\documentclass{article}
\usepackage[utf8]{inputenc}

\pdfoutput=1
\usepackage{mathtools}
\usepackage{amsfonts}
\usepackage{amssymb}
\usepackage{centernot}
\usepackage{graphicx}
\usepackage{amsthm}
\usepackage{enumerate}

\usepackage{tikz}
\usepackage{tikz-cd}
\usetikzlibrary{decorations.markings}
\usetikzlibrary{arrows}
\usetikzlibrary{calc}
\providecommand{\U}[1]{\protect\rule{.1in}{.1in}}
\newtheorem{theorem}{Theorem}

\newtheorem{corollary}[theorem]{Corollary}

\newtheorem{definition}[theorem]{Definition}

\newtheorem{lemma}[theorem]{Lemma}

\newtheorem{proposition}[theorem]{Proposition}

\theoremstyle{remark}

\newcommand{\Mr}{M_A^{\textup{red}}}
\newcommand{\Me}{M_A^{\textup{enr}}}

\newcommand{\ann}{\textup{ann}}

\allowdisplaybreaks

\begin{document}

\title{Peripheral elements in reduced Alexander modules}
\author{Lorenzo Traldi\\Lafayette College\\Easton, PA 18042, USA\\traldil@lafayette.edu
}
\date{ }
\maketitle

\begin{abstract}
We discuss meridians and longitudes in reduced Alexander modules of classical and virtual links. When these elements are suitably defined, each link component will have many meridians, but only one longitude. Enhancing the reduced Alexander module by singling out these peripheral elements provides a significantly stronger link invariant. In particular, the enhanced module determines all linking numbers in a link; in contrast, the module alone does not even detect how many linking numbers are $0$.


\end{abstract}

\section{Introduction}

We begin by establishing notation and terminology.

We use the term \emph{link diagram} to refer to an oriented, virtual link diagram, i.e., a subset $D$ of $\mathbb R ^2$ obtained in the following way. Begin with a finite number of oriented, piecewise smooth closed curves, $C_1, \dots, C_ \mu$, with only finitely many (self-) intersections, all of which are crossings (transverse double points). Designate each crossing as either \emph{classical} or \emph{virtual}. At each classical crossing, remove a short piece of the underpassing segment on each side of the crossing. Then draw a small circle around each virtual crossing. A link diagram is \emph{classical} if all of its crossings are classical.

Two link diagrams are \emph{equivalent} if they are related through a finite sequence of four types of changes: the three classical Reidemeister moves, which affect only classical crossings, and the detour move, in which a segment of a curve $C_i$ with no classical crossing is replaced by any other piecewise smooth segment with the same endpoints, also with no classical crossing. All four kinds of moves preserve the integrity of the images of $C_1, \dots, C_ \mu$, so it is reasonable to say that an equivalence class of link diagrams represents a \emph{link} $L=K_1 \cup \dots \cup K_ \mu$, where each component $K_i$ is represented in a diagram by the closed curve $C_i$ with the same index. (N.b.\ Our links are oriented, virtual links.)  It is well known that equivalence classes of classical link diagrams represent actual links in $\mathbb S^3$ \cite{R}, and equivalence classes of general link diagrams represent actual links in thickened surfaces \cite{K}, but in this paper we think of links simply as equivalence classes of link diagrams. We refer to the literature (e.g., the text of Manturov and Ilyutko \cite{MI}) for background on links and their diagrams.

If $D$ is a link diagram, then the result of removing the short pieces of underpassing segments near classical crossings is to cut the original curves into arcs. N.b.\ Only the underpassing segment is cut at a classical crossing; the overpassing segment is not cut. Also, the segments involved in a virtual crossing simply pass straight through the crossing; neither segment is cut, and the two segments are not considered to be attached to each other at the crossing. The set of arcs in $D$ is denoted $A(D)$, and the set of classical crossings in $D$ is denoted $C(D)$. (We do not adopt notation for virtual crossings because as Kauffman said \cite{vkt}, for our purposes they are ``not really there.'') Let $\kappa_D:A(D) \to \{1, \dots, \mu \}$ be the map with $\kappa_D(a)=i$ if and only if $a$ is part of the image of $K_i$ in $D$.

\begin{figure} 
\centering
\begin{tikzpicture} [>=angle 90]
\draw [thick] (1,.5) -- (0.8,0.4);
\draw [thick] [<-] (0.8,0.4) -- (-1,-.5);
\draw [thick] (-1,.5) -- (-.2,0.1);
\draw [thick] (0.2,-0.1) -- (1,-.5);
\node at (1.5,0.5) {$a(c)$};
\node at (-1.5,0.5) {$b_2(c)$};
\node at (1.5,-0.5) {$b_1(c)$};
\end{tikzpicture}
\caption{A classical crossing $c$.}
\label{crossfig}
\end{figure}
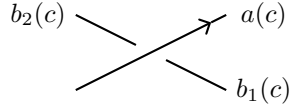

Now, let $\Lambda= \mathbb Z [t,t^{-1}]$ be the ring of Laurent polynomials in $t$, with integer coefficients. If $D$ is a diagram of $L=K_1 \cup \dots \cup K_ \mu$, let $\Lambda^{A(D)}$ and $\Lambda^{C(D)}$ be the free $\Lambda$-modules on the sets $A(D)$ and $C(D)$, respectively. Let $a,b_1,b_2:C(D) \to A(D)$ be the functions indicated in Fig.\ \ref{crossfig}. (Notice that $b_1(c)$ is on the right-hand side of $a(c)$, and $b_2(c)$ is on the left-hand side of $a(c)$.) Then there is a $\Lambda$-linear map $\varrho_D:\Lambda^{C(D)} \to \Lambda^{A(D)}$ given by
\[
\varrho_D(c)=(1-t)a(c)+tb_1(c)-b_2(c) \quad \forall c \in C(D).
\]
 The cokernel of $\varrho_D$ is the \emph{reduced Alexander module} of $D$, denoted $\Mr(D)$. There is a canonical map onto the quotient $\varsigma_D:\Lambda^{A(D)} \to \Mr(D)$, with $\ker \varsigma_D = \varrho_D(\Lambda^{C(D)})$.

Let $\epsilon: \Lambda \to \mathbb Z$ be the augmentation map, i.e., the ring homomorphism with $\epsilon (t^{\pm 1}) = 1$. Also, let $\mathbb Z _ \epsilon$ be the $\Lambda$-module structure on the abelian group $\mathbb Z$ determined by $\epsilon$, so that $\lambda \cdot n = \epsilon(\lambda) n \thickspace \allowbreak \forall \lambda \in \Lambda \thickspace \allowbreak \forall n \in \mathbb Z$. 

The definition of $\varrho_D$ given above implies that there is a $\Lambda$-linear map $\varphi_D: \Mr(D) \to \Lambda \oplus (\mathbb Z_\epsilon) ^{\mu -1}$, defined as follows. If $\kappa_D(a)=1$, then $\varphi_D(\varsigma_D(a))=(1,0,\dots,0)$; and if $\kappa_D(a)=i>1$, then $\varphi_D(\varsigma_D(a))=(1,0,\dots,0, 1, 0, \dots, 0)$, with the second $1$ in the $i$th coordinate.

We take a moment to briefly discuss some of the history of these ideas. Many researchers have studied Alexander matrices of classical link diagrams since they were introduced almost 100 years ago \cite{A}. For the first few decades, research involving Alexander matrices seems to have been focused on the associated Alexander polynomials and elementary ideals; see \cite{F} for instance. The fact that Alexander matrices describe modules over Laurent polynomial rings does not seem to have been studied systematically before Crowell's work \cite{C1, C2a, C2, C3}.  

A few years later, other researchers studied peripheral structures within multivariate Alexander modules. For instance, Levine's work on two-component link modules \cite{L} included longitudinal elements, although he only defined them when the linking number is $0$. A thorough presentation of classical multivariate Alexander module theory is given in Hillman's book \cite{H}. A notable aspect of the theory is that substantial results involving the longitudes seem to require the hypothesis that the linking numbers are all $0$, or the even stronger hypothesis that the $(\mu -1)$st elementary ideal is $0$; see \cite[Sec.\ 4.7]{H}.

The results of the present paper indicate that in contrast, the peripheral elements of reduced Alexander modules have useful properties no matter what the linking numbers are. In particular, our longitudinal elements are given by an explicit formula (Definition \ref{longdef}), they are indeterminate (each component of $L$ has a unique longitude in $\Mr(L)$), and they play a distinctive structural role in the module $\Mr(L)$ (see Section \ref{torsion}).

Here is our definition of meridianal elements in $\Mr(L)$.

\begin{definition}
\label{merdef}
An element $x \in \Mr(D)$ is a \emph{meridian of $K_i$ in $\Mr(D)$} if $\varphi_D(x)=\varphi_D(\varsigma_D(a))$ for any arc $a \in A(D)$ with $\kappa_D(a)=i$. The set of meridians of $K_i$ in $\Mr(D)$ is denoted $M_i(D)$.
\end{definition}

Notice that the map $\varphi_D$ is crucial in Definition \ref{merdef}. As far as we know, this map is not part of the classical theory discussed above; $\varphi_D$ was derived from the multivariate version of the link module sequence of Crowell \cite{C1} in recent work regarding the connection between Alexander modules and quandles. In particular, a special case of Definition \ref{longdef} in \cite{mvaq4} indicated that a certain kind of link diagram gives rise to ``longitudes'' in the reduced Alexander module, which define special automorphisms of the medial quandle: displacements that act as identity maps on individual quandle orbits. The present paper developed when we wondered about the significance of the ideas of \cite{mvaq4} outside the theory of link quandles. As this is our motivation, we do not focus on quandles here.

\begin{figure} 
\centering
\begin{tikzpicture} [>=angle 90]
\draw [thick] (1,.5) -- (0.8,0.4);
\draw [thick] [<-] (0.8,0.4) -- (-1,-.5);
\draw [thick] (-1,.5) -- (-.2,0.1);
\draw [thick] [<-] (0.8,-0.4) -- (0.2,-0.1);
\draw [thick] (0.8,-0.4) -- (1,-.5);
\draw [thick] (5,.5) -- (4.8,0.4);
\draw [thick] [<-] (4.8,0.4) -- (3,-.5);
\draw [thick]  [<-] (3.2,.4) -- (3.8,0.1);
\draw [thick]  (3,0.5) -- (3.8,0.1);
\draw [thick] (5,-.5) -- (4.2,-0.1);
\node at (0,-1) {$w(c)=-1$};
\node at (4,-1) {$w(c)=1$};
\end{tikzpicture}
\caption{The writhe of a classical crossing $c$.}
\label{writhefig}
\end{figure}
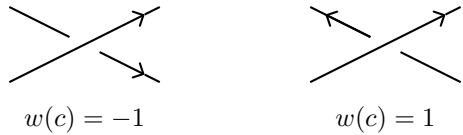

In order to define longitudes in $\Mr(D)$ we use a familiar notion, the writhe of a classical crossing; see Fig.\ \ref{writhefig}. 

\begin{definition}
\label{longdef}
If $1 \leq i \leq \mu$, then the \emph{longitude of $K_i$ in $\Mr(D)$} is given by the formula
\[
\chi_i(D) =  \varsigma_D \left (\sum_{\substack{c \in C(D)\\ \kappa_D(b_1(c))=i}} w(c)a(c) - \frac{1}{2} \sum_{\substack{c \in C(D)\\ \kappa_D(b_1(c))=i}} w(c)(b_1(c)+b_2(c)) \right).
\]
\end{definition}

Here are three comments on Definition \ref{longdef}.

1. The coefficient $\frac{1}{2}$ requires explanation, because $\Lambda$ has no such element. Each $a \in A(D)$ occurs either $0$ times or $2$ times as $b_1(c)$ or $b_2(c)$ for some classical crossing $c$ with $\kappa_D(b_1(c))=i$. Therefore 
\[
\sum_{\substack{c \in C(D)\\ \kappa_D(b_1(c))=i}} w(c)(b_1(c)+b_2(c)) = \sum_{a \in A(D)} m_a a \text{,}
\]
where each coefficient $m_a$ is $-2$, $0$, or $2$. These coefficients can be multiplied by $\frac{1}{2}$ in the obvious ways, so $\chi_i(D)$ is well-defined.

2. A longitude of a classical link $L=K_1 \cup \dots \cup K_ \mu$ is usually defined to have linking number $0$ with the corresponding component $K_i$. Definition \ref{longdef} tells us that $\chi_i$ instead represents a curve whose linking number with the entire link $L$ is $0$, in the sense that
\[
\chi_i(D) = \varsigma_D \left(\sum_{a \in A(D)} n_a a \right)
\]
for some integers $n_a$ with $\sum n_a=0$. The reader might prefer to think of $\chi_i$ as a ``totally unlinked parallel'' rather than a ``longitude.''

3. We use the symbol $\chi_i$ simply because symbols related to the letter \emph{l} are preempted: we use $\lambda$ for elements of $\Lambda$, and $\ell$ for linking numbers.

\begin{definition}
\label{enh}
The \emph{enhanced} reduced Alexander module of a link diagram $D$ is the list
\[
\Me(D)=(\Mr(D),M_1(D), \dots, M_\mu(D),\chi_1(D), \dots , \chi_\mu(D)).
\]
\end{definition}

Notice that for each link component, $\Me(D)$ includes a set of meridians  but only one longitude. This is a fundamental difference between Definition \ref{enh} and the familiar theory of peripheral elements in classical link groups, where the meridian-longitude pairs are defined up to simultaneous conjugation.

\begin{definition}
\label{enhiso}
Let $D$ and $D'$ be $\mu$-component link diagrams, and suppose $f:\Mr(D) \to \Mr(D')$ is an isomorphism of $\Lambda$-modules. If $f(M_i(D))=M_i(D')$ and $f(\chi_i(D))=\chi_i(D')$ for every $i \in \{ 1, \dots, \mu \}$, then we say $f$ is an isomorphism between $\Me(D)$ and $\Me(D')$.
\end{definition}

There is another way to state the requirement $f(M_i(D)) = M_i(D')$ in Definition \ref{enhiso}. The sets $M_1(D), \dots, M_\mu(D)$ generate the $\Lambda$-module $\Mr(D)$, and the map $\varphi_D$ is constant on each set $M_i(D)$, so once $M_1(D), \dots, M_\mu(D)$ are determined, the map $\varphi_D$ is determined too. It follows that an isomorphism $f:\Mr(D) \to \Mr(D')$ has $f(M_i(D))=M_i(D') \thickspace \allowbreak \forall i \in \{ 1, \dots, \mu \}$ if and only if it has $\varphi_D = \varphi_{D'} \circ f$.

Our primary result is that up to isomorphism, enhanced reduced Alexander modules are link type invariants:
\begin{theorem}
\label{main}
Let $D$ and $D'$ be equivalent link diagrams. Then there is an isomorphism $f:\Me(D) \to \Me(D')$.
\end{theorem}

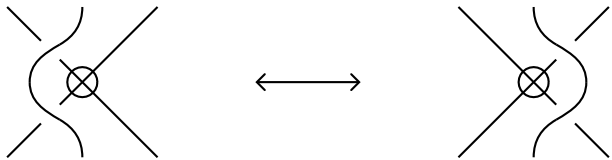
\begin{figure} 
\centering
\begin{tikzpicture} [>=angle 90]
\draw [thick] [<-] (-.7,0) -- (0,0);
\draw [thick] [->] (0,0) -- (.7,0);
\draw [thick] (-3.3,-.3) -- (-2,1);
\draw [thick] (-3.55,-.55) -- (-4,-1);
\draw [thick] (-3.3,.3) -- (-2,-1);
\draw [thick] (-3.55,.55) -- (-4,1);
\draw [thick] (-3,1) to [out=-90, in = 30] (-3.3,.5);
\draw [thick] (-3.3,.5) to [out=-150, in = 90] (-3.7,0);
\draw [thick] (-3,-1) to [out=90, in = -30] (-3.3,-.5);
\draw [thick] (-3.7,0) to [out=-90, in = 150] (-3.3,-.5);
\draw  [thick, domain=0:360] plot ({-3+.2*cos(\x)}, {.2*sin(\x)});
\draw [thick] (3.3,-.3) -- (2,1);
\draw [thick] (3.55,-.55) -- (4,-1);
\draw [thick] (3.3,.3) -- (2,-1);
\draw [thick] (3.55,.55) -- (4,1);
\draw [thick] (3,1) to [out=-90, in = 150] (3.3,.5);
\draw [thick] (3.3,.5) to [out=-30, in = 90] (3.7,0);
\draw [thick] (3,-1) to [out=90, in = -150] (3.3,-.5);
\draw [thick] (3.7,0) to [out=-90, in = 30] (3.3,-.5);
\draw  [thick, domain=0:360] plot ({3+.2*cos(\x)}, {.2*sin(\x)});
\end{tikzpicture}
\caption{The welded move.}
\label{welded}
\end{figure}

We are grateful to an anonymous reader for the observation that $\Me(D)$ is also invariant under the welded overpass move of \cite{FRR}, pictured in Fig. \ref{welded}.

After proving Theorem \ref{main} we often drop $D$ from our notation, and write $\Me(L)=(\Mr(L),M_1(L), \dots, M_\mu(L),\chi_1(L), \dots , \chi_\mu(L))$. 

It turns out that the enhanced reduced Alexander module is a much more sensitive link invariant than the reduced Alexander module alone. One instance of this increased sensitivity involves linking numbers. The classical definition of linking numbers was extended to virtuals by Goussarov, Polyak and Viro \cite{GPV} as follows.

\begin{definition}
\label{lno}
Let $D$ be a diagram of $L=K_1 \cup \dots \cup K_{\mu}$. If $i \neq j \in \{1, \dots, \mu \}$, then the \emph{linking number} of $K_j$ over $K_i$ is
\[
\ell_{j/i}(K_i,K_j) = \ell_{j/i}(K_j,K_i)=\sum_{\substack{c \in C(D)\\ \kappa_D(a(c))=j, \kappa_D(b_1(c))=i}} w(c).
\]
\end{definition}

Notice that we use the subscript $j/i$ to reflect the fact that the linking number $\ell_{j/i}(K_i,K_j) = \ell_{j/i}(K_j,K_i)$ is determined by the crossings of $K_j$ over $K_i$. In the classical case this is not important, because the linking number is symmetric in $i$ and $j$; but for virtuals, the linking number of $K_j$ over $K_i$ is independent from the linking number of $K_i$ over $K_j$. Different authors have used different notation to reflect this lack of symmetry; for instance Goussarov, Polyak and Viro \cite{GPV} used $lk_{j/i}$ rather than $\ell_{j/i}(K_i,K_j)= \ell_{j/i}(K_j,K_i)$, and Chrisman \cite {Ch} used $v \ell k(K_j,K_i)$.

It is not hard to find examples of links that have isomorphic reduced Alexander modules despite having different arrays of linking numbers. In contrast, enhanced reduced Alexander modules have the following property.

\begin{theorem}
\label{linkingthm}
If $L=K_1 \cup \dots \cup K_ \mu$ and $L'=K'_1 \cup \dots \cup K'_ {\mu'}$ are links with $\Me(L) \cong \Me(L')$, then $\mu=\mu'$ and the linking numbers in $L$ are precisely the same as the linking numbers in $L'$.
\end{theorem}

The sensitivity of $\Me(L)$ as a link invariant does not lie only in the fact that it determines linking numbers. For instance, many link invariants fail to distinguish the two versions of the Borromean rings -- they have isomorphic Alexander modules, link module sequences, groups and fundamental quandles -- but it is not difficult to show that $\Me$ distinguishes them. See Section \ref{examples} for details.

Here is an outline of our discussion. We prove Theorem \ref{main} in Section \ref{Reidemeister}, and Theorem \ref{linkingthm} in Section \ref{lnos}. In Section \ref{torsion}, we point out that the longitudes have a special role in the module $\Mr(L)$: they generate the submodule annihilated by $1-t$. In Section \ref{knots} we observe that when $L$ is a knot, enhancing $\Mr(L)$ with peripheral information is of little value. In Section \ref{examples}, we mention some examples. In Section \ref{questions}, we ask some questions about possible extensions of our results.

\section{Proof of Theorem \ref{main}}
\label{Reidemeister}

In this section we prove Theorem \ref{main}, by analyzing the effects of classical Reidemeister moves on $\Me(D)$. It is not necessary to discuss detour moves or welded moves, because they do not affect any of the information used to define $\Me(D)$. We give details for four particular instances of Reidemeister moves. Polyak \cite{P} showed that all other Reidemeister moves can be obtained from these four.

\begin{figure} [bth]
\centering
\begin{tikzpicture} [>=angle 90]
\draw [thick] [->] (-0.25,-0.25) -- (.5,.5);
\draw [thick] (.5,.5) -- (.75,.75);
\draw [thick] (-0.25,-0.25) to [out=225, in=0] (-1,-.75);
\draw [thick] (-1,-.75) to [out=180, in=-90] (-1.9,0);
\draw [thick] (-1,.75) to [out=180, in=90] (-1.9,0);
\draw [thick] (-1,.75) to [out=0, in=-45] (-0.2,0.2);
\draw [thick] (.75,-0.75) -- (0.1,-0.1);
\draw [thick] [->] (-5,-0.75) -- (-5,.5);
\draw [thick] (-5,.75) -- (-5,.5);
\node at (-5.4,0) {$a_0$};
\node at (1.4,-0.7) {$b_1(c_0)$};
\node at (2.4,0.6) {$a_0=a(c_0)=b_2(c_0)$};
\node at (-.5,0) {$c_0$};
\end{tikzpicture}
\caption{An $\Omega.1$ move changes $D$ into $D'$.}
\label{firstmove}
\end{figure}
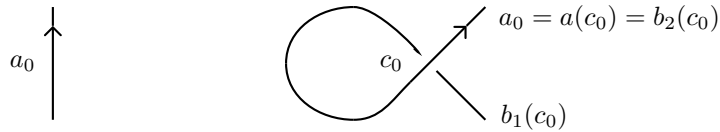

Suppose $D'$ is obtained from $D$ by performing the $\Omega.1$ move illustrated in Fig.\ \ref{firstmove}. Let $F:A(D') \to A(D)$ be the function with $F(a_0)=F(b_1(c_0))=a_0$ and $F(a)=a$ for every other arc. Then $F$ extends to a $\Lambda$-linear surjection $\Lambda^{A(D')} \to \Lambda^{A(D)}$, which we also denote $F$.

Notice that $F(\varrho_{D'}(c_0))=F(-ta_0+tb_1(c_0))=0$, and $F(\varrho_{D'}(c)) = \varrho_D(c)$ for every other crossing. As $b_1(c_0)-a_0=t^{-1}\varrho_{D'}(c_0) \in \varrho_{D'}(\Lambda^{C(D')})$ , it follows that $F$ induces an isomorphism $f:\Mr(D') \to \Mr(D)$, with $f \varsigma_{D'} = \varsigma_D F$. This equality implies $f(M_i(D'))=M_i(D) \thickspace \allowbreak \forall i \in \{ 1, \dots, \mu \}$. Moreover the only difference between the defining formulas of the longitudes of $D$ and $D'$ is that if $i=\kappa_D(a_0)$ then $\chi_i(D')$ includes a contribution from $c_0$. This difference is insignificant because the contribution from $c_0$ is
\[
\varsigma_{D'}(w(c_0)a(c_0) - \frac{1}{2} w(c_0)(b_1(c_0)+b_2(c_0))) 
\]
\[
= w(c_0) \varsigma_{D'}(a_0 - \frac{1}{2} (a_0+a_0)) = w(c_0)\varsigma_{D'}(0)=0.
\]

\begin{figure} 
\centering
\begin{tikzpicture} [>=angle 90]
\draw [thick]  (-0.25,-0.25) -- (.3,.3);
\draw [thick] [->] (.75,.75) -- (.3,.3);
\draw [thick] (-0.25,-0.25) to [out=225, in=0] (-1,-.75);
\draw [thick] (-1,-.75) to [out=180, in=-90] (-1.9,0);
\draw [thick] (-1,.75) to [out=180, in=90] (-1.9,0);
\draw [thick] (-1,.75) to [out=0, in=-45] (-0.2,0.2);
\draw [thick] (.75,-0.75) -- (0.1,-0.1);
\draw [thick] (-5,-0.75) -- (-5,.3);
\draw [thick] [->] (-5,.75) -- (-5,.3);
\node at (-5.4,0) {$a_0$};
\node at (1.4,-0.7) {$b_2(c_0)$};
\node at (2.4,0.6) {$a_0=a(c_0)=b_1(c_0)$};
\node at (-.5,0) {$c_0$};
\end{tikzpicture}
\caption{An $\Omega.1$ move changes $D$ into $D'$.}
\label{firstmoveb}
\end{figure}
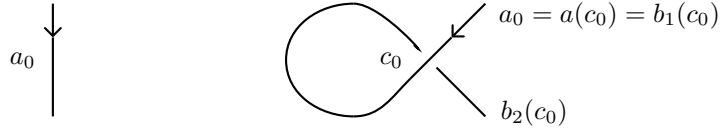

Now, suppose $D'$ is obtained from $D$ by performing the $\Omega.1$ move illustrated in Fig.\ \ref{firstmoveb}. Let $F:A(D') \to A(D)$ be the function with $F(a_0)=F(b_2(c_0))=a_0$ and $F(a)=a$ for every other arc. As before, $F$ extends to a $\Lambda$-linear surjection $F:\Lambda^{A(D')} \to \Lambda^{A(D)}$, which has $F(\varrho_{D'}(c_0))=F(a_0-b_2(c_0))=0$ and $F(\varrho_{D'}(c)) = \varrho_D(c)$ for every other crossing. The equality $a_0-b_2(c_0)=\varrho_{D'}(c_0)$ implies that just as before, $F$ induces an isomorphism $f:\Mr(D') \to \Mr(D)$ with $f \varsigma_{D'} = \varsigma_D F$. As before, this equality implies $f(M_i(D'))=M_i(D) \thickspace \allowbreak \forall i \in \{ 1, \dots, \mu \}$. Again, the only difference between the defining formulas of the longitudes of $D$ and $D'$ is the contribution from $c_0$ to $\chi_{\kappa_D(a_0)}(D')$, and again, this difference is insignificant:
\[
\varsigma_{D'}(w(c_0)a(c_0) - \frac{1}{2} w(c_0)(b_1(c_0)+b_2(c_0))) 
\]
\[
= w(c_0) \varsigma_{D'}(a_0 - \frac{1}{2} (a_0+a_0)) = w(c_0)\varsigma_{D'}(0)=0.
\]

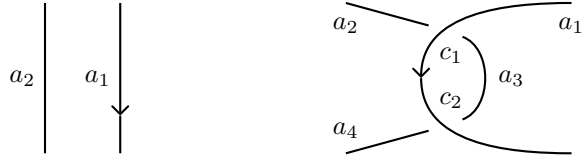
\begin{figure} 
\centering
\begin{tikzpicture} [>=angle 90]
\draw [thick] (1,-1) to [out=-180, in=-90] (-1,0);
\draw [thick] [->] (1,1) to [out=-180, in=90] (-1,0);
\draw [thick] (-2,-1) -- (-.9,-0.7);
\draw [thick] (-2,1) -- (-.9,0.7);
\draw [thick] (-0.45,-0.55) to [out=20, in=-20] (-0.45,0.55);
\draw [thick] (-6,-1) -- (-6,1);
\draw [thick] [<-] (-5,-.5) -- (-5,1);
\draw [thick] (-5,-.5) -- (-5,-1);
\node at (-0.6,0.3) {$c_1$};
\node at (-0.6,-0.3) {$c_2$};
\node at (1,0.7) {$a_1$};
\node at (-2,0.7) {$a_2$};
\node at (.2,0) {$a_3$};
\node at (-2,-0.7) {$a_4$};
\node at (-6.3,0) {$a_2$};
\node at (-5.3,0) {$a_1$};
\end{tikzpicture}
\caption{An $\Omega.2$ move changes $D$ into $D'$.}
\label{secondmove}
\end{figure}

Now, suppose $D'$ is obtained from $D$ by performing an $\Omega.2$ move, as in Fig.\ \ref{secondmove}. Let  $F:\Lambda^{A(D')} \to \Lambda^{A(D)}$ be the $\Lambda$-linear surjection with $F(a_3)=(1-t)a_1+ta_2$, $F(a_4)=a_2$ and $F(a)=a$ otherwise. Notice that $F( \varrho_{D'} (c_1)) = F( \varrho_{D'} (c_2))=(1-t)a_1 +ta_2 - ((1-t)a_1+ta_2)=0$, and $F(\varrho_{D'}(c)) = \varrho_D(c)$
for every other crossing $c$, so $F$ induces a $\Lambda$-linear map $f:\Mr(D') \to \Mr(D)$ with $f \varsigma_{D'} = \varsigma_D f$. The image of $F$ includes $\varsigma_D(a)$ for every $a \in A(D)$, so $f$ is surjective. In addition, $f$ is injective, because 
\[
\varsigma_{D'}(a_2) =\varsigma_{D'}(a_2+0) =  \varsigma_{D'}(a_2+t^{-1}\varrho_{D'}(c_2-c_1))= \varsigma_{D'}(a_2+a_4-a_2)=\varsigma_{D'}(a_4).
\]
Thus $f$ is an isomorphism.

The fact that $\varphi_D \varsigma_D F(a_3) = \varphi_D \varsigma_D F(a_4) =\varphi_D(a_2)$ implies $f(M_i(D')) = M_i(D) \thickspace \allowbreak \forall i \in \{ 1, \dots, \mu \}$. The only difference between the defining formulas of the longitudes of $D$ and $D'$ is that if $i=\kappa_D(a_2)$, then $\chi_i(D')$ includes contributions from $c_1$ and $c_2$. This difference is insignificant because the sum of these contributions is
\[
\varsigma_{D'}(w(c_1)a(c_1) + w(c_2)a(c_2)- \frac{1}{2} (w(c_1)(b_1(c_1)+b_2(c_1))+w(c_2)(b_1(c_2)+b_2(c_2))) 
\]
\[
=\varsigma_{D'}(w(c_1)a_1 - w(c_1)a_1- \frac{1}{2} (w(c_1)(a_2+a_3)-w(c_1)(a_4+a_3)))
\]
\[
=\varsigma_{D'}(w(c_1) \frac{1}{2} (a_4-a_2)) \text{,}
\]
which is $0$ because $\varsigma_{D'}(a_2) = \varsigma_{D'}(a_4)$.

\begin{figure} [b]
\centering
\begin{tikzpicture} [>=angle 90]  
\draw [thick] [->] (2.2,-1.2) -- (1.2,-.7); 
\draw [thick] (1.2,-.7) -- (-2.6,1.2); 
\draw [thick] (-.08,.06) -- (2.2,1.2); 
\draw [thick] [->]  (-.32,-.06) -- (-1.6,-.7); 
\draw [thick] (-2.6,-1.2) --(-1.6,-.7);
\draw [thick] (-2.6,0) to [out=0, in=210] (-1.75,.63);
\draw [thick] (-1.52,.78) to [out=30,in=150] (1.12,.78);
\draw [thick] (2.2,0) to [out=180, in=-30] (1.35,.63);
\draw [thick] [->] (-4.8,-1.2) -- (-8.6,.7);
\draw [thick] (-8.6,.7) -- (-9.6,1.2);
\draw [thick] [<-] (-5.8,.7) -- (-4.8,1.2);
\draw [thick]  (-7.08,.06) -- (-5.8,.7);
\draw [thick] (-7.32,-.06) -- (-9.6,-1.2);
\draw [thick] (-9.6,0) to [out=0, in=-210] (-8.75,-.63);
\draw [thick] (-8.52,-.78) to [out=-30,in=-150] (-5.88,-.78);
\draw [thick] (-4.8,0) to [out=180, in=30] (-5.65,-.63);
\node at (-4.98,-.84) {$a_1$};
\node at (-5.4,0) {$a_6$};
\node at (-4.98,.78) {$a_2$};
\node at (-7.2,-.9) {$a_5$};
\node at (-9.48,-.84) {$a_3$};
\node at (-9,0) {$a_4$};
\node at (2.02,-.84) {$a_1$};
\node at (1.6,0) {$a_6$};
\node at (2.02,.78) {$a_2$};
\node at (-.2,.9) {$a_7$};
\node at (-2.48,-.84) {$a_3$};
\node at (-2,0) {$a_4$};
\node at (-.2,-.3) {$c$};
\node at (-7.2,0.3) {$c$};
\end{tikzpicture}
\caption{An $\Omega.3$ move changes $D$ into $D'$.}
\label{thirdmove2}
\end{figure}
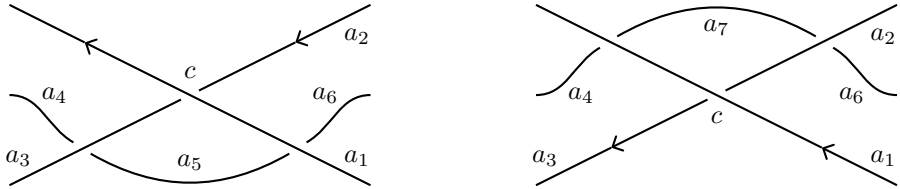

Now, let $D$ and $D'$ be the diagrams on the left and right of Fig.\ \ref{thirdmove2}. The images under $\varrho_D$ of the generators of $\Lambda^{C(D)}$ corresponding to the three pictured crossings of $D$ are $\varrho_1=(1-t)a_1+ta_2-a_3$, $\varrho_2=(1-t)a_1+ta_6-a_5$, and $\varrho_3=(1-t)a_3+ta_4-a_5$.
The images under $\varrho_{D'}$ of the generators of $\Lambda^{C(D')}$ corresponding to the three pictured crossings of $D'$ are $\varrho'_1=(1-t)a_1+ta_2-a_3$, $\varrho'_2=(1-t)a_1+ta_7-a_4$, and $\varrho'_3=(1-t)a_2+ta_7-a_6$.

Let $F:\Lambda^{A(D)} \to \Lambda^{A(D')}$ be the $\Lambda$-linear map with $F(a)=a$ for $a \neq a_5$, and  $F(a_5)=a_4+t(a_6-a_7)$. Then $F$ is surjective, because its image includes every $a \neq a_7 \in A(D')$, and also includes $a_7=F(t^{-1} \cdot(a_4-a_5+ta_6))$. $F$ is injective, too: if $x \neq 0 \in \Lambda^{A(D)}$ then either the $a_5$ coordinate of $x$ is nonzero, in which case the $a_7$ coordinate of $F(x)$ is nonzero, or else the $a_5$ coordinate of $x$ is $0$, in which case $x$ and $F(x)$ are precisely the same as linear combinations of generators. We conclude that $F$ is an isomorphism of $\Lambda$-modules. 

Notice that $F(\varrho_1)=\varrho'_1$, $F(\varrho_2)=\varrho'_2$, and $F(\varrho_3)=(t-1)\varrho'_1+(1-t)\varrho'_2+t\varrho'_3$, so $
F(\varrho_D(\Lambda^{C(D)}))=\varrho_{D'}(\Lambda^{C(D')})$. It follows that $F$ induces an isomorphism $f:\Mr(D) \to \Mr(D')$. As $\varphi_{D'}f\varsigma_D(a_5) = \varphi_{D'}\varsigma_{D'}(a_4+t(a_6-a_7))= \varphi_{D'}\varsigma_{D'}(a_4)$, $f(M_i(D))=M_i(D') \thickspace \allowbreak \forall i \in \{ 1, \dots, \mu \}$.


The only difference between the defining formulas of longitudes of $D$ and $D'$ in Fig.\ \ref{thirdmove2} involves the pictured crossings not marked $c$. Suppose the crossing furthest to the left in $D$ has writhe $w$. If $i=\kappa_D(a_4)$, then the contribution of the two unmarked crossings to $\chi_i(D)$ is 
\[
w \cdot \varsigma_D(a_3-a_1 - \frac{1}{2}(a_4-a_6)).
\]
The contribution of the two unmarked crossings to $\chi_i(D')$ is
\[
w \cdot \varsigma_{D'}(a_2-a_1 - \frac{1}{2}(a_6-a_4)).
\]
The isomorphism $f$ matches these two contributions precisely:
\[
f(w \cdot \varsigma_D(a_3-a_1 - \frac{1}{2}(a_4-a_6))) = w \cdot \varsigma_{D'}(a_3-a_1 - \frac{1}{2}(a_4-a_6)))
\]
\[
=w \cdot \varsigma_{D'}(a_2-a_1 - \frac{1}{2}(a_6-a_4)) + w \cdot \varsigma_{D'}(a_3-a_2+a_6-a_4)
\]
\[
=w \cdot \varsigma_{D'}(a_2-a_1 - \frac{1}{2}(a_6-a_4)) + w \cdot \varsigma_{D'}(\varrho'_2-\varrho'_1-\varrho'_3) \text{,}
\]
and of course $\varsigma_{D'}(\varrho'_2-\varrho'_1-\varrho'_3)=0$.

\section{Linking numbers}
\label{lnos}
Recall the map $\varphi_L: \Mr(L) \to \Lambda \oplus (\mathbb Z_\epsilon) ^{\mu -1}$ mentioned in the introduction. Let $\epsilon^\mu: \Lambda \oplus (\mathbb Z_\epsilon) ^{\mu -1} \to \mathbb Z ^ \mu$ be the map given by $\epsilon:\Lambda \to \mathbb Z$ in the first coordinate, and identity maps in the later coordinates.
 
Definitions \ref{longdef} and \ref{lno} immediately imply the following.
 
\begin{lemma}
\label{linklem}
Let $L=K_1 \cup \dots \cup K_{\mu}$ be a link. If $i \in \{1, \dots, \mu\}$, then $\epsilon^\mu \varphi_L (\chi_i(L))$ is the element of $\mathbb Z ^ \mu$ whose $j$th coordinate is $\ell_{j/i}(K_i,K_j)$ if $i \neq j$, and 
\[
-\sum _{k \neq i} \ell_{k/i}(K_i,K_k)
\]
if $i=j$.
\end{lemma}

Theorem \ref{linkingthm} follows readily from Lemma \ref{linklem}. Suppose $L=K_1 \cup \dots \cup K_{\mu}$ and $L'=K'_1 \cup \dots \cup K'_{\mu}$ are links, and  $f:\Me(L) \to \Me(L')$ is an isomorphism. As noted after Definition \ref{enhiso}, it follows that $\varphi_L = \varphi_{L'} f$, and hence $\epsilon^\mu \varphi_L = \epsilon^\mu \varphi_{L'} f$. Lemma \ref{linklem} then implies that $\ell_{j/i}(K_i,K_j)=\ell_{j/i}(K'_i,K'_j) \thickspace \allowbreak \forall i \neq j \in \{ 1, \dots, \mu \}$. This completes the proof of Theorem \ref{linkingthm}.

\section{Longitudes and torsion}
\label{torsion}

The following result was mentioned in \cite{mvaq4}. The discussion there involved quandles and a special case of Definition \ref{longdef}, which required a special type of link diagram. The proof given below applies to any diagram, and does not involve quandles.

\begin{theorem}
\label{fed}
The submodule of $\Mr(L)$ generated by $\chi_1(L), \dots, \chi_\mu(L)$ is $\ann(1-t) = \{ x \in M_A^{\textup{red}}(L) \mid (1-t)x=0\}$.
\end{theorem}
\begin{proof} Let $D$ be a diagram of $L$. 

We begin by verifying that $(1-t)\chi_i(D)$ is always $0$. According to the description of $\Mr(D)$ given in the introduction, every $c \in C(D)$ has
\[
\varsigma_D((1-t)a(c)) = \varsigma_D(b_2(c)-tb_1(c)).
\]
There are elements $\lambda_a \in \Lambda$, for $a \in A(D)$, such that 
\[
(1-t)\chi_i(D) =  (1-t)\varsigma_D \left (\sum_{\substack{c \in C(D)\\ \kappa_D(b_1(c))=i}} w(c)a(c) - \frac{1}{2} \sum_{\substack{c \in C(D)\\ \kappa_D(b_1(c))=i}} w(c)(b_1(c)+b_2(c)) \right)
\]
\[
=\varsigma_D \left( \sum_{\substack{c \in C(D)\\ \kappa_D(b_1(c))=i}} w(c) \cdot  (b_2(c)-tb_1(c)) - (1-t) \cdot \frac{1}{2}\sum_{\substack{c \in C(D)\\ \kappa_D(b_1(c))=i}} w(c)(b_1(c)+b_2(c))\right)
\]
\[
=\varsigma_D \left(\sum_{a \in A(D)} \lambda_a a \right).
\]

Consider an arc $a \in A(D)$. If $a$ is not an underpassing arc at any crossing, then of course $\lambda_a=0$. Otherwise, let $c_1,c_2$ be the classical crossings with $a \in \{b_1(c_i),b_2(c_i)\}$, with $a$ oriented from $c_1$ toward $c_2$. (It is possible that $c_1=c_2$.) As illustrated in Fig.\ \ref{bifig}, there are four possible configurations for the arc $a$ and the crossings $c_1,c_2$. We check that each of the four configurations leads to $\lambda_a=0$. 
\begin{figure} 
\centering
\begin{tikzpicture} [>=angle 90]
\draw [thick] [->] (-6,-.5) -- (-6,.5);
\draw [thick] [->] (-5.9,0) -- (-5,0);
\draw [thick] (-6.3,0) -- (-6.1,0);
\draw [thick] (-5,0) -- (-4.1,0);
\draw [thick] [->] (-4,.5) -- (-4,-.5);
\draw [thick] (-3.7,0) -- (-3.9,0);
\node at (-5.6,-.3) {$c_1$};
\node at (-4.4,-.3) {$c_2$};
\draw [thick] (-3.3,0) -- (-3.1,0);
\draw [thick] [->] (-3,-.5) -- (-3,.5);
\draw [thick] [->] (-2.9,0) -- (-2,0);
\draw [thick] (-2,0) -- (-1.1,0);
\draw [thick] [->] (-1,-.5) -- (-1,.5);
\draw [thick] (-.7,0) -- (-.9,0);
\node at (-2.6,-.3) {$c_1$};
\node at (-1.4,-.3) {$c_2$};
\draw [thick] (-.1,0) -- (-.3,0);
\draw [thick] [->] (0,.5) -- (0,-.5);
\draw [thick] [->] (0.1,0) -- (1,0);
\draw [thick] (1,0) -- (1.9,0);
\draw [thick] [->] (2,.5) -- (2,-.5);
\draw [thick] (2.3,0) -- (2.1,0);
\node at (0.4,-.3) {$c_1$};
\node at (1.6,-.3) {$c_2$};
\draw [thick] (2.7,0) -- (2.9,0);
\draw [thick] [->] (3,.5) -- (3,-.5);
\draw [thick] [->] (3.1,0) -- (4,0);
\draw [thick] (4,0) -- (4.9,0);
\draw [thick] [->] (5,-.5) -- (5,.5);
\draw [thick] (5.3,0) -- (5.1,0);
\node at (3.4,-.3) {$c_1$};
\node at (4.6,-.3) {$c_2$};
\node at (4,.3) {$a$};
\node at (1,.3) {$a$};
\node at (-2,.3) {$a$};
\node at (-5,.3) {$a$};
\end{tikzpicture}
\caption{The four configurations of classical crossings at the ends of an arc $a$.}
\label{bifig}
\end{figure}
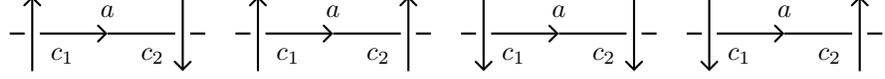

In the first configuration $a=b_1(c_1)=b_1(c_2)$ and $w(c_1)=-1=-w(c_2)$, so 
\[
\lambda_a=(--t-t)-(1-t) \frac{1}{2}(-1+1)=0-0=0.
\]
In the second configuration $a=b_1(c_1)=b_2(c_2)$ and $w(c_1)=-1=w(c_2)$, so 
\[
\lambda_a=(--t-1)-(1-t) \frac{1}{2}(-1-1)=(t-1)+(1-t)=0.
\]
In the third configuration $a=b_2(c_1)=b_1(c_2)$ and $w(c_1)=1=w(c_2)$, so 
\[
\lambda_a=(1-t)-(1-t) \frac{1}{2}(1+1)=(1-t)-(1-t)=0.
\]
In the fourth configuration $a=b_2(c_1)=b_2(c_2)$ and $w(c_1)=1=-w(c_2)$, so 
\[
\lambda_a=(1-1)-(1-t) \frac{1}{2}(1-1)=0-0=0.
\]
It follows that $(1-t)\chi_i(D)=\varsigma_D(\sum_{a \in A(D)} \lambda_a a)=0$, as required.

Now, suppose $x \in \Mr(L)$ has $(1-t)x=0$. We must show that $x$ is equal to a linear combination of $\chi_1(D), \dots, \chi_\mu(D)$.

Let $x' \in \Lambda^{A(D)}$ have $x= \varsigma_D(x')$.
Then $\varsigma_D((1-t)x')=(1-t)x=0$, so $(1-t)x' \in \ker \varsigma_D$; hence there is a function $g_{x'}:C(D) \to \Lambda$ such that
\begin{equation}
\label{lastprop1}
(1-t)x' = \varrho_D \left( \sum_{c \in C(D)} g_{x'}(c) \cdot c \right).
\end{equation}

We claim that $x'$ can be chosen so that $g_{x'}(c)$ is always an integer. To verify the claim, note first that if $c \in C(D)$ then $\epsilon(\epsilon(g_{x'}(c))-g_{x'}(c))=0$, so there is an element $\lambda_c \in \Lambda$ such that $\epsilon(g_{x'}(c))-g_{x'}(c)=\lambda_c (1-t)$. Let 
\[
x''=x'+\varrho_D \left( \sum_{c \in C(D)} \lambda_c \cdot c \right).
\]
Notice that
\[
(1-t)x'' = 
\varrho_D \left(\sum_{c \in C(D)} (g_{x'}(c) + \lambda_c(1-t) ) \cdot c \right) = \varrho_D \left(\sum_{c \in C(D)} \epsilon(g_{x'}(c)) \cdot c \right) \text{,}
\]
so $g_{x''}(c) = \epsilon(g_{x'}(c))$ is always an integer. As $\varsigma_D(x'')=\varsigma_D(x')=x$, it follows that the claim is satisfied if we choose $x''$ to play the role of $x'$.

Recalling (\ref{lastprop1}), we now have 
\[
(1-t)x' = \sum_{c \in C(D)} g_{x'}(c) \cdot \varrho_D (c) = \sum_{c \in C(D)} g_{x'}(c) \cdot ((1-t)a(c)+tb_1(c)-b_2(c))
\]
\begin{equation}
\label{lastprop2}
=(1-t)\sum_{c \in C(D)} g_{x'}(c) a(c)+ \sum_{c \in C(D)} g_{x'}(c) \cdot (tb_1(c)-b_2(c)).
\end{equation}
It follows that the last sum displayed above is a multiple of $1-t$:
\[
\sum_{c \in C(D)} g_{x'}(c) \cdot (tb_1(c)-b_2(c)) = (1-t) \cdot \left( x' - \sum_{c \in C(D)} g_{x'}(c) a(c) \right)
\]
\[
=(1-t) \cdot  \sum_{a \in A(D)} f(a) a
\]
for some function $f:A(D) \to \Lambda$. This equality holds in the free $\Lambda$-module $\Lambda^{A(D)}$, so the sums must match precisely. That is, for each $a \in A(D)$ this equality must hold: 
\begin{equation}
\label{lastprop3}
t \cdot \sum_{\substack{c \in C(D)\\ b_1(c)=a}} g_{x'}(c) - \sum_{\substack{c \in C(D)\\ b_2(c)=a}} g_{x'}(c) = (1-t)f(a).
\end{equation}

Suppose $a \in A(D)$; then there cannot be more than two crossings with $a \in \{b_1(c),b_2(c)\}$. Recall that $g_{x'}(c) \in \mathbb Z \thickspace \allowbreak \forall c \in C(D)$.  Clearly then for (\ref{lastprop3}) to hold, the left-hand side of (\ref{lastprop3}) must be an integer multiple of $1-t$. Therefore one of the following must be true: (a) there are crossings $c,c' \in C(D)$ with $g_{x'}(c)=g_{x'}(c')=-f(a)$ and $b_1(c)=a=b_2(c')$, or (b) there are crossings $c,c' \in C(D)$ with $g_{x'}(c)=-g_{x'}(c')$, $f(a)=0$ and either $b_1(c)=a=b_1(c')$ or $b_2(c)=a=b_2(c')$. Case (a) includes the possibility that $c=c'$.

Consulting Fig.\ \ref{bifig}, we see that $w(c)=w(c')$ in case (a), and $w(c)=-w(c')$ in case (b). In either case, we have $g_{x'}(c)w(c) = g_{x'}(c')w(c')$. This equality can also be stated as follows: for each arc $a \in A(D)$, there is an integer $m_a$ such that every classical crossing at which $a$ is an underpassing arc has $g_{x'}(c)w(c) = m_a$. Notice that if $a$ and $a'$ are underpassing arcs at the same classical crossing $c$, then $m_a=g_{x'}(c)w(c)=m_{a'}$. Walking from crossing to crossing along the arcs of $D$, we deduce that the value of $m_a$ is constant on each component $K_i$ of $L$. Denote this constant value $m_i$. 

As $w(c)$ is always $\pm 1$, the equality $g_{x'}(c)w(c) = m_a$ tells us that $g_{x'}(c) = m_a w(c)$. In case (a), it follows that
\[
f(a) = -g_{x'}(c)=-g_{x'}(c') = -\frac{1}{2}(g_{x'}(c)+g_{x'}(c'))= - m_a \cdot \frac{1}{2}(w(c)+w(c')).
\]
In case (b) we have $g_{x'}(c)=-g_{x'}(c')$, so
\[
f(a) = 0 = -\frac{1}{2}(g_{x'}(c)+g_{x'}(c')) = - m_a \cdot \frac{1}{2}(w(c)+w(c')).
\]

Notice that we have the same formula for $f(a)$ in either case. Using this formula and $g_{x'}(c) = m_a w(c)$ to rewrite equation (\ref{lastprop2}), we obtain the following.
\[
(1-t)x'=(1-t) \cdot \sum_{c \in C(D)} g_{x'}(c) a(c)+ \sum_{c \in C(D)} g_{x'}(c) \cdot (tb_1(c)-b_2(c))
\]
\[
=(1-t) \cdot \sum_{c \in C(D)} g_{x'}(c) a(c)+ (1-t) \cdot \sum_{a \in A(D)} f(a) a
\]
\[
=(1-t) \cdot \sum_{i=1}^ \mu \sum_{\substack{c \in C(D)\\ \kappa_D(b_1(c))=i}} m_i w(c) a(c)+ (1-t) \cdot \sum_{a \in A(D)} - m_a \cdot \frac{1}{2}(w(c)+w(c')) a
\]
\[
=(1-t) \cdot \sum_{i=1}^ \mu m_i \sum_{\substack{c \in C(D)\\ \kappa_D(b_1(c))=i}} w(c) a(c)- (1-t) \cdot \sum_{i=1}^ \mu m_i \cdot \frac{1}{2} \sum_{\substack{c \in C(D)\\ \kappa_D(b_1(c))=i}} w(c) (b_1(c)+b_2(c))
\]
\begin{equation}
\label{lastprop4}
=(1-t) \cdot \sum_{i=1}^ \mu m_i \cdot \left( \sum_{\substack{c \in C(D)\\ \kappa_D(b_1(c))=i}} w(c) a(c)- \frac{1}{2} \sum_{\substack{c \in C(D)\\ \kappa_D(b_1(c))=i}} w(c) (b_1(c)+b_2(c)) \right).
\end{equation}

As the equality (\ref{lastprop4}) holds in the free $\Lambda$-module $\Lambda^{A(D)}$, we can cancel the $1-t$ factors to obtain 
\[
x' =\sum_{i=1}^ \mu m_i \cdot \left( \sum_{\substack{c \in C(D)\\ \kappa_D(b_1(c))=i}} w(c) a(c)- \frac{1}{2} \sum_{\substack{c \in C(D)\\ \kappa_D(b_1(c))=i}} w(c) (b_1(c)+b_2(c)) \right).
\]
It follows that $x=\varsigma_D(x') = \sum m_i \chi_i(D)$. \end{proof}

\section{Knots}
\label{knots}
In this section we focus on knots, i.e., links with $\mu=1$. It turns out that enhancing the reduced Alexander module of a knot $K$ is of no value; $\Me(K)$ is completely determined by the $\Lambda$-module $\ker \varphi_K$, the \emph{Alexander invariant} of $K$. To prove this, we need two special properties of the Alexander invariants of knots: multiplication by $1-t$ always defines an automorphism, and longitudinal elements are always equal to $0$. These properties are well known in the classical case \cite[ps.\ 6 and 34]{H}. Proofs are included for the sake of completeness.

\begin{lemma}
\label{knotlem}
If $K$ is a knot, then $\ker \varphi_K = (1-t) \cdot \ker \varphi_K$.
\end{lemma}
\begin{proof}
Let $D$ be a diagram of $K$. We index the arcs of $D$ as $a_1, \dots, a_n$, in such a way that we encounter the arcs in order as we walk along $K$ in $D$. 

Let $x \in \ker \varphi_K$. Then 
\begin{equation}
\label{knoteq}
x=\varsigma_D \left( \sum_{i=1}^n \lambda_i a_i \right)
\end{equation}
for some elements $\lambda_1, \dots, \lambda_n \in \Lambda$ with $\sum \lambda_i=\varphi_K(x) = 0$.

Recall that $\epsilon:\Lambda \to \mathbb Z$ is the map with $\epsilon(t ^{\pm 1})=1$; then $\ker \epsilon$ is the principal ideal of $\Lambda$ generated by $1-t$. 

Now, suppose it is not possible to find $\lambda_1, \dots, \lambda_n \in \Lambda$ such that equation (\ref{knoteq}) holds and $\epsilon(\lambda_i)=0 \thickspace \allowbreak \forall i \in \{ 1, \dots, n \}$. Choose $\lambda_1, \dots, \lambda_n$ so that the least index $i_0$ with $\epsilon(\lambda_{i_0}) \neq 0$ is as large as possible.

Case 1: Suppose $i_0<n$. Then there is a crossing $c$ of $D$ with $\{b_1(c),b_2(c)\}=\{a_{i_0},a_{1+i_0}\}$. Let $a(c)=a_j$.

If $b_1(c)=a_{i_0}$ and $b_2(c)=a_{1+i_0}$, then we have 
\[
0=-\epsilon(\lambda_{i_0}) \cdot 0 =  -\epsilon(\lambda_{i_0}) \cdot \varsigma_D(\varrho_D(c))=-\epsilon(\lambda_{i_0}) \cdot \varsigma_D((1-t)a(c)+tb_1(c)-b_2(c))
\]
\[
=-\epsilon(\lambda_{i_0})  (1-t) \cdot \varsigma_D(a_j) -\epsilon(\lambda_{i_0})t \cdot \varsigma_D(a_{i_0})+\epsilon(\lambda_{i_0}) \cdot \varsigma_D(a_{1+i_0}).
\]
It follows that if we add $-\epsilon(\lambda_{i_0})  (1-t)$ to $\lambda_j$, add $-\epsilon(\lambda_{i_0})t$ to $\lambda_{i_0}$, and add $\epsilon(\lambda_{i_0})$  to $\lambda_{1+i_0}$, then equation (\ref{knoteq}) still holds. This contradicts our choice of $\lambda_1, \dots, \lambda_n$, because $\epsilon(\lambda_1), \dots, \epsilon(\lambda_{i_0})$ are now all $0$.

It follows by contradiction that $b_1(c)=a_{1+i_0}$ and $b_2(c)=a_{i_0}$. Then we have 
\[
0=\epsilon(\lambda_{i_0}) \cdot 0 =  \epsilon(\lambda_{i_0}) \cdot \varsigma_D(\varrho_D(c))=\epsilon(\lambda_{i_0}) \cdot \varsigma_D((1-t)a(c)+tb_1(c)-b_2(c))
\]
\[
=\epsilon(\lambda_{i_0})  (1-t) \cdot \varsigma_D(a_j) +\epsilon(\lambda_{i_0})t \cdot \varsigma_D(a_{1+i_0})-\epsilon(\lambda_{i_0}) \cdot \varsigma_D(a_{i_0}).
\]
Therefore if we add $\epsilon(\lambda_{i_0})  (1-t)$ to $\lambda_j$, add $\epsilon(\lambda_{i_0})t$ to $\lambda_{1+i_0}$, and add $-\epsilon(\lambda_{i_0})$  to $\lambda_{i_0}$, then equation (\ref{knoteq}) still holds. This again contradicts our choice of $\lambda_1, \dots, \lambda_n$, because $\epsilon(\lambda_1), \dots, \epsilon(\lambda_{i_0})$ are now all $0$.

Case 2: Suppose $i_0=n$. That is, $\epsilon(\lambda_1), \dots, \epsilon(\lambda_{n-1})$ are all $0$, and $\epsilon(\lambda_n) \neq 0$. This is impossible, as $\sum \lambda_i = 0$.

We conclude by contradiction that it is possible to choose $\lambda_1, \dots, \lambda_n$ so that equation (\ref{knoteq}) holds and $\epsilon(\lambda_i)=0 \thickspace \allowbreak \forall i \in \{ 1, \dots, n \}$. Each $\lambda_i$ is a multiple of $1-t$; say $\lambda_i=(1-t)\lambda'_i$. Then $0= \sum \lambda_i = (1-t)\sum \lambda'_i$. As $\Lambda$ is an integral domain, it follows that $0= \sum \lambda'_i$, so 
\[
x' = \varsigma_D \left( \sum_{i=1}^n \lambda'_i a_i \right)
\]
is an element of $\ker \varphi_K$ with $x=(1-t)x'$.
\end{proof}

\begin{corollary}
\label{knotcor}
If $K$ is a knot, then multiplication by $1-t$ defines an automorphism of $\ker \varphi_K$. Also, $\chi_1(K)=0$.
\end{corollary}
\begin{proof}
Multiplication by $1-t$ defines an endomorphism of any $\Lambda$-module, of course. Lemma \ref{knotlem} tells us that for $\ker \varphi_K$, the endomorphism is surjective. The ring $\Lambda$ is Noetherian, so the finitely generated $\Lambda$-module $\ker \varphi_K$ is Noetherian too. Therefore every surjective endomorphism of $\ker \varphi_K$ is an automorphism.

Theorem \ref{fed} tells us that $(1-t)\chi_1(K)=0$, so it follows that $\chi_1(K)=0$.
\end{proof}

\begin{theorem}
\label{knotthm}
Let $K$ and $K'$ be knots. Then $\Me(K) \cong \Me(K')$ if and only if $\ker \varphi_K \cong \ker \varphi_{K'}$.
\end{theorem}

\begin{proof}
Let $D$ be a diagram of $K$. The elements of $M_1(K)$ include every $\varsigma_D(a)$ with $a \in A(D)$, and these elements generate $\Mr(K)$, so $\ker \varphi_K$ is the submodule of $\Mr(K)$ generated by $\{x-y \mid x,y \in M_1(K) \}$. The same holds for $K'$.

If $\Me(K) \cong \Me(K')$, then there is an isomorphism $f:\Mr(K) \to \Mr(K')$ of $\Lambda$-modules, with $f(M_1(K))=M_1(K')$. It follows that 
\[
f(\{x-y \mid x,y \in M_1(K) \})= \{x-y \mid x,y \in M_1(K') \}
\]
and hence according to the preceding paragraph, $f(\ker \varphi_K) = \ker \varphi_{K'}$.

For the converse, observe that the $\Lambda$-linear epimorphism $\varphi_K:\Mr(K) \to \Lambda$ must split. That is, there is an isomorphism $g:\Mr(K) \to \Lambda \oplus \ker \varphi_K$ such that $\varphi_K$ is the composition of $g$ and the projection map $\Lambda \oplus \ker \varphi_K \to \Lambda$. There is a similar isomorphism $g':\Mr(K) \to \Lambda \oplus \ker \varphi_{K'}$.

Now, suppose there is an isomorphism $\ker \varphi_K \to \ker \varphi_{K'}$. It defines an isomorphism $f:\Lambda \oplus \ker \varphi_K \to \Lambda \oplus \ker \varphi_{K'}$, which is the identity map on the first coordinate. It follows that $f$ maps $g(M_1(K))$, the subset of $\Lambda \oplus \ker \varphi_K$ containing all ordered pairs of the form $(1,x)$, to $g'(M_1(K'))$. Therefore $h=(g')^{-1}fg:\Mr(K) \to \Mr(K')$ is an isomorphism with $h(M_1(K))=M_1(K')$. 

Corollary \ref{knotcor} implies that $h(\chi_1(K))=h(0)=0=\chi_1(K')$, so $h$ defines an isomorphism between $\Me(K)$ and $\Me(K')$.
\end{proof}
\section{Examples}
\label{examples}

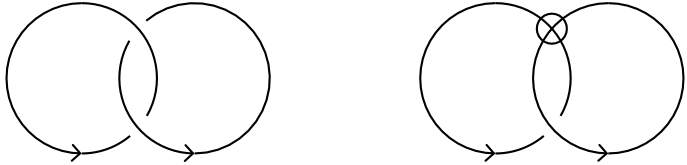
\begin{figure} [bth]
\centering
\begin{tikzpicture} [>=angle 90]
\draw [->] [thick, domain=-270:-90] plot ({-3.5+cos(\x)}, {sin(\x)});
\draw [thick, domain=-50:-90] plot ({-3.5+cos(\x)}, {sin(\x)});
\draw [thick, domain=-30:90] plot ({-3.5+cos(\x)}, {sin(\x)});
\draw [->] [thick, domain=-210:-90] plot ({-2+cos(\x)}, {sin(\x)});
\draw  [thick, domain=-90:130] plot ({-2+cos(\x)}, {sin(\x)});
\draw [->] [thick, domain=-270:-90] plot ({3.5+cos(\x)}, {sin(\x)});
\draw [thick, domain=-90:90] plot ({3.5+cos(\x)}, {sin(\x)});
\draw [->] [thick, domain=-270:-90] plot ({2+cos(\x)}, {sin(\x)});
\draw  [thick, domain=-90:-50] plot ({2+cos(\x)}, {sin(\x)});
\draw  [thick, domain=-30:90] plot ({2+cos(\x)}, {sin(\x)});
\draw  [thick, domain=0:360] plot ({2.75+.2*cos(\x)}, {.66+.2*sin(\x)});
\end{tikzpicture}
\caption{Classical and virtual Hopf links.}
\label{hopfs}
\end{figure}

In the abstract, we made the comment that in general, $\Mr(L)$ does not detect how many linking numbers in $L$ are $0$. This comment is easy to verify for virtual links. For instance, both links pictured in Fig.\ \ref{hopfs} have $\Mr(L) \cong \Lambda \oplus (\Lambda/(1-t))$, even though the classical Hopf link has $\ell_{1/2}(K_1,K_2)=\ell_{2/1}(K_1,K_2)=1$ and the virtual Hopf link has $\ell_{1/2}(K_1,K_2)=0$ and $\ell_{2/1}(K_1,K_2)=1$. The fact that classical links have $\ell_{j/i}(K_i,K_j)=\ell_{i/j}(K_i,K_j)$ implies that classical examples of the comment must have $\mu >2$. Two appropriate examples were discussed in \cite[Sec.\ 4.2]{mvaq4}; they are 3-component classical links with isomorphic reduced Alexander modules, one with linking numbers $1,0,0$ and the other with linking numbers $2,2,-1$.

\begin{figure} 
\centering
\begin{tikzpicture} [>=angle 90]
\draw [->] [thick, domain=-240:-90] plot ({-4+1.5*cos(\x)}, {-1.4+1.5*sin(\x)});
\draw [thick, domain=-90:-10] plot ({-4+1.5*cos(\x)}, {-1.4+1.5*sin(\x)});
\draw [thick, domain=10:100] plot ({-4+1.5*cos(\x)}, {-1.4+1.5*sin(\x)});
\draw [->] [thick, domain=10:180] plot ({-5+1.5*cos(\x)}, {1.5*sin(\x)});
\draw [thick, domain=180:240] plot ({-5+1.5*cos(\x)}, {1.5*sin(\x)});
\draw [thick, domain=260:350] plot ({-5+1.5*cos(\x)}, {1.5*sin(\x)});
\draw [->] [thick, domain=236:360] plot ({-3+1.5*cos(\x)}, {1.5*sin(\x)});
\draw [thick, domain=140:220] plot ({-3+1.5*cos(\x)}, {1.5*sin(\x)});
\draw [thick, domain=0:124] plot ({-3+1.5*cos(\x)}, {1.5*sin(\x)});
\node at (-6.5,-2.7) {$B$};
\node at (-4.7,-1.7) {$v$};
\node at (-6.2,.5) {$u$};
\node at (-4.7,.5) {$x$};
\node at (-1.9,.5) {$w$};
\node at (-2.6,-.3) {$z$};
\node at (-4.5,-2.5) {$y$};
\node at (-6.4,1.2) {$K_1$};
\node at (-1.6,1.2) {$K_2$};
\node at (-2.7,-2.7) {$K_3$};
\draw [thick, domain=-240:-90] plot ({2.5+1.5*cos(\x)}, {-1.4+1.5*sin(\x)});
\draw [<-] [thick, domain=-90:-10] plot ({2.5+1.5*cos(\x)}, {-1.4+1.5*sin(\x)});
\draw [thick, domain=10:100] plot ({2.5+1.5*cos(\x)}, {-1.4+1.5*sin(\x)});
\draw [thick, domain=10:180] plot ({1.5+1.5*cos(\x)}, {1.5*sin(\x)});
\draw [<-] [thick, domain=180:240] plot ({1.5+1.5*cos(\x)}, {1.5*sin(\x)});
\draw [thick, domain=260:350] plot ({1.5+1.5*cos(\x)}, {1.5*sin(\x)});
\draw [thick, domain=236:360] plot ({3.5+1.5*cos(\x)}, {1.5*sin(\x)});
\draw [thick, domain=140:220] plot ({3.5+1.5*cos(\x)}, {1.5*sin(\x)});
\draw [<-] [thick, domain=0:124] plot ({3.5+1.5*cos(\x)}, {1.5*sin(\x)});
\node at (5,-2.7) {$B'$};
\node at (1.8,-1.7) {$u'$};
\node at (.3,.5) {$v'$};
\node at (1.8,.5) {$w'$};
\node at (4.6,.5) {$x'$};
\node at (3.9,-.3) {$y'$};
\node at (2,-2.5) {$z'$};
\node at (0.1,1.2) {$K'_1$};
\node at (4.9,1.2) {$K'_2$};
\node at (1.2,-2.7) {$K'_3$};
\end{tikzpicture}
\caption{$B$ and $B'$ are the two versions of the Borromean rings.}
\label{borr}
\end{figure}
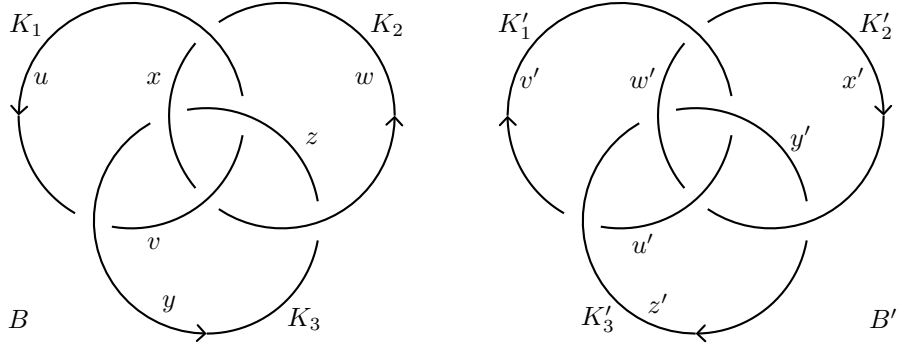

In the introduction, we mentioned that in general, $\Me(L)$ is more sensitive than the combination of $\Mr(L)$ and linking numbers. In the rest of this section, we illustrate this point with a well-known example, the Borromean rings, pictured in Fig.\ \ref{borr}. It is a good exercise to verify that the Borromean rings are equivalent to their mirror image, and also equivalent to a link obtained by reversing orientations of two components. These observations bring the number of possibly distinct oriented versions of the link down to the two pictured ones, but they do not verify that these two are distinct. N.b.\ As mentioned in the introduction, our notion of ``equivalent'' requires that the indices of link components be preserved. If we were to allow permutations of the component indices, then there would be only one link -- to obtain $B'$ from $B$, rotate the diagram of $B$ through an angle of $\pi$ around a vertical axis. This rotation interchanges the indices of the first two link components, of course.

As a small abuse of notation, we use $B$ and $B'$ to denote both the link diagrams pictured in Fig.\ \ref{borr} and the links represented by the diagrams. 

Notice that if we tabulate information regarding arcs and crossings in $B$ and $B'$, we get precisely the same data, except for the use of apostrophes in $B'$. For instance, there is only one crossing $c$ of $B$ with $a(c)=u$, and it has $b_1(c)=w$ and $b_2(c)=x$; there is also only one crossing $c'$ of $B'$ with $a(c')=u'$, and it has $b_1(c')=w'$ and $b_2(c')=x'$. Furthermore, the two diagrams have the same association between link components and arcs; for instance, the arcs corresponding to $K_2$ are $w$ and $x$, and the arcs corresponding to $K'_2$ are $w'$ and $x'$. 

Therefore, it is impossible to distinguish $B$ from $B'$ using any link invariant determined by the functions $a,b_1,b_2$ mapping crossings to arcs, along with the function $\kappa$ mapping arcs to $\{1, \dots, \mu \}$. There are many very powerful such invariants, including not only $\Mr$ but also the link group and the fundamental quandle. Nevertheless, we can distinguish $B$ from $B'$ by showing that $\Me(B) \not \cong \Me(B')$. Here are the details.

The $\Lambda$-module $\Mr(B)$ is generated by six elements, $\varsigma_B(u)$ through $\varsigma_B(z)$, subject to relations corresponding to the crossings in $B$. Here are the crossing relations:
\[
\varsigma_B(v) = (1-t)\varsigma_B(y) + t\varsigma_B(u),\varsigma_B(x) = (1-t)\varsigma_B(u) + t\varsigma_B(w),
\]
\[
\varsigma_B(z) = (1-t)\varsigma_B(w) + t\varsigma_B(y),\varsigma_B(v) = (1-t)\varsigma_B(z) + t\varsigma_B(u),
\]
\[
\varsigma_B(x) = (1-t)\varsigma_B(v) + t\varsigma_B(w),\varsigma_B(z) = (1-t)\varsigma_B(x) + t\varsigma_B(y)
\]
We use the first three crossing relations to eliminate the generators $\varsigma_B(v)$, $\varsigma_B(x)$ and $\varsigma_B(z)$. It follows that $\Mr(B)$ is generated by $\varsigma_B(u)$, $\varsigma_B(w)$ and $\varsigma_B(y)$, subject to relations that result from the last three crossing relations:
\[
(1-t)\varsigma_B(y) + t\varsigma_B(u) = (1-t)((1-t)\varsigma_B(w) + t\varsigma_B(y)) + t\varsigma_B(u)
\]
\[
(1-t)\varsigma_B(u) + t\varsigma_B(w) = (1-t)((1-t)\varsigma_B(y) + t\varsigma_B(u)) + t\varsigma_B(w)
\]
\[
(1-t)\varsigma_B(w) + t\varsigma_B(y) = (1-t)((1-t)\varsigma_B(u) + t\varsigma_B(w)) + t\varsigma_B(y)
\]

These relations are equivalent to $(1-t)^2\varsigma_B(y-w) = 0$, $(1-t)^2\varsigma_B(y-u) = 0$ and $(1-t)^2\varsigma_B(w-u) = 0$ (respectively). Therefore
\begin{equation}
\label{modiso}
\Mr(B) \cong \Lambda \oplus (\Lambda/(1-t)^2) \oplus (\Lambda/(1-t)^2) \text{,}
\end{equation}
with the three summands generated by $\varsigma_B(u)$, $\varsigma_B(w-u)$ and $\varsigma_B(y-u)$, respectively. We deduce the following.

\begin{proposition}
\label{annprops}
The annihilator of $1-t$ in $\Mr(B)$ is
\[
\ann(1-t)=\{(1-t) \cdot (\lambda_1\varsigma_B(w-u)+\lambda_2\varsigma_B(y-u)) \mid \lambda_1,\lambda_2 \in \Lambda\}.
\]
Moreover, if $\lambda_1,\lambda_2,\widetilde \lambda_1,\widetilde \lambda_2 \in \Lambda$  then 
\[
(1-t) \cdot (\lambda_1\varsigma_B(w-u)+\lambda_2\varsigma_B(y-u))=(1-t) \cdot (\widetilde \lambda_1\varsigma_B(w-u)+\widetilde \lambda_2\varsigma_B(y-u))
\]
if and only if $\epsilon(\lambda_1)=\epsilon(\widetilde \lambda_1)$ and $\epsilon(\lambda_2)=\epsilon(\widetilde \lambda_2)$.
\end{proposition}
\begin{proof}
In the direct sum on the right-hand side of (\ref{modiso}), it is obvious that $\ann(1-t)$ is the set of 3-tuples of the form $(0,(1-t)\lambda_1 + (1-t)^2,(1-t)\lambda_2 + (1-t)^2)$, where $\lambda_1,\lambda_2 \in \Lambda$. As the sum is direct, two such 3-tuples are equal if and only if their coordinates are equal. If $\lambda,\widetilde \lambda \in \Lambda$ then the cosets $(1-t)\lambda + (1-t)^2$ and $(1-t) \widetilde \lambda + (1-t)^2$ are equal if and only if $\lambda-\widetilde \lambda$ is a multiple of $1-t$, and this is true if and only if $\epsilon(\lambda)=\epsilon(\widetilde \lambda)$. 

The statement of the proposition is the result of translating the observations of the preceding paragraph into $\Mr(B)$ via the isomorphism (\ref{modiso}).
\end{proof}

Let $\pi_1:\Mr(B) \to \Lambda$ be the first coordinate of the isomorphism (\ref{modiso}). The isomorphism (\ref{modiso}) tells us that every $m \in \Mr(B)$ may be written in the form $m=\pi_1(m) \varsigma_B(u) + \lambda_1 \varsigma_B(w-u) + \lambda_2 \varsigma_B(y-u)$ for some $\lambda_1,\lambda_2 \in \Lambda$, and this expression is unique except for the fact that any multiple of $(1-t)^2$ may be added to $\lambda_1$ or $\lambda_2$. We then have the following formula for $\varphi_B:\Mr(B) \to \Lambda \oplus (\mathbb Z _ \epsilon)^2$.
\[
\varphi_B(m)=\varphi_B(\pi_1(m) \varsigma_B(u) +\lambda_1 \varsigma_B(w-u) + \lambda_2 \varsigma_B(y-u)) 
\]
\begin{equation}
\label{phiform}
= \pi_1(m)(1,0,0) + \lambda_1(0,1,0)+\lambda_2(0,0,1)= (\pi_1(m),\epsilon(\lambda_1),\epsilon(\lambda_2))
\end{equation}
We can now restate Proposition \ref{annprops} in the following way.

\begin{corollary}
\label{anncor}
The annihilator of $1-t$ in $\Mr(B)$ is
\[
\ann(1-t)=\{(1-t) m \mid m \in \Mr(B) \text{ and } \pi_1(m) = 0 \}.
\]
Moreover, if $m,\widetilde m \in \Mr(B)$ and $\pi_1 (m)=0=\pi_1 (\widetilde m)$, then $(1-t)m=(1-t) \widetilde m$ if and only if $\varphi_B(m)= \varphi_B(\widetilde m)$.
\end{corollary}
\begin{proof}
The isomorphism (\ref{modiso}) tells us that $\pi_1(m)=0$ if and only if $m$ is of the form $m=\lambda_1 \varsigma_B(w-u) + \lambda_2 \varsigma_B(y-u)$ for some $\lambda_1,\lambda_2 \in \Lambda$. This fact allows us to deduce the equality $\ann(1-t)=\{(1-t) m \mid m \in \Mr(B) \text{ and } \pi_1(m) = 0 \}$ from Proposition \ref{annprops}.

Now, suppose $m,\widetilde m \in \Mr(B)$ and $\pi_1 (m)=0=\pi_1 (\widetilde m)$. Then there are $\lambda_1,\lambda_2,\widetilde \lambda_1,\widetilde \lambda_2 \in \Lambda$ with $m=\lambda_1 \varsigma_B(w-u) + \lambda_2 \varsigma_B(y-u)$ and $\widetilde m=\widetilde \lambda_1 \varsigma_B(w-u) + \widetilde \lambda_2 \varsigma_B(y-u)$. Formula (\ref{phiform}) tells us that $\varphi_B(m)=\varphi_B(\widetilde m)$ if and only if the equalities $\epsilon(\lambda_1)=\epsilon(\widetilde \lambda_1)$ and $\epsilon(\lambda_2)=\epsilon(\widetilde \lambda_2)$ are both true, and Proposition \ref{annprops} tells us that these equalities are both true if and only if $(1-t)m= (1-t)\widetilde m$.
\end{proof}

The discussion above also applies to $\Mr(B')$, with apostrophes attached to symbols associated with $B'$. For instance, $\pi'_1:\Mr(B') \to \Lambda$ is the first coordinate of the $B'$ version of the isomorphism (\ref{modiso}).

We claim that there is no isomorphism $f:\Mr(B) \to \Mr(B')$ such that $f(M_1(B)) = M_1(B')$, $f(M_2(B)) = M_2(B')$, $f(M_3(B)) = M_3(B')$ and $f(\chi_1(B)) = \chi_1(B')$. To verify this claim, note first that
\begin{equation}
\label{longform1}
\chi_1(B)=\varsigma_B(y)-\varsigma_B(z) = \varsigma_B(y)-(1-t)\varsigma_B(w) - t\varsigma_B(y)
=(1-t)m \text{,}
\end{equation}
where $m=(-1)\cdot\varsigma_B(w-u)+ 1 \cdot\varsigma_B(y-u)$. Similarly,
\begin{equation}
\label{longform2}
\chi_1(B')=\varsigma_{B'}(z')-\varsigma_B(y') = (1-t) \varsigma_{B'}(w') + t\varsigma_{B'}(y')-\varsigma_{B'}(y')=(1-t)m' \text{,}
\end{equation}
where $m'=1 \cdot\varsigma_{B'}(w'-u')+ (-1) \cdot \varsigma_B(y'-u')$.


Now, suppose $f:\Mr(B) \to \Mr(B')$ is an isomorphism with $f(M_1(B)) = M_1(B')$, $f(M_2(B)) = M_2(B')$ and $f(M_3(B)) = M_3(B')$. As noted in the introduction, it follows that $\varphi_B=\varphi_{B'}f$. According to formula (\ref{phiform}), $\varphi_B=\varphi_{B'}f$ implies that $\pi_1 = \pi'_1 f$. Therefore, $f(m)$ and $m'$ are elements of $\Mr(B')$ with $\pi'_1(f(m))=\pi_1(m)=0=\pi'_1(m')$ and
\[
\varphi_{B'}(f(m))=\varphi_B(m) = (0,-1,1) \neq (0,1,-1) = \varphi_{B'}(m').
\]
According to the $B'$ version of Corollary \ref{anncor}, the inequality $\varphi_{B'}(f(m)) \neq \varphi_{B'}(m')$ implies
that 
\[
f(\chi_1(B)) =f((1-t)m) = (1-t)f(m)\neq (1-t)m'=\chi_1(B').
\]
This verifies our claim.

Notice by the way that there is certainly an isomorphism $f:\Mr(B) \to \Mr(B')$ with $f(M_1(B)) = M_1(B')$, $f(M_2(B)) = M_2(B')$ and $f(M_3(B)) = M_3(B')$. Such an $f$ may be defined simply by attaching apostrophes to symbols: $f(\varsigma_B(u))=\varsigma_{B'}(u')$, etc. Also, it turns out that then $-f:\Mr(B) \to \Mr(B')$ is an isomorphism with $-f(\chi_1(B)) = \chi_1(B')$, $-f(\chi_2(B)) = \chi_2(B')$ and $-f(\chi_3(B)) = \chi_3(B')$. But as we have seen, there is no isomorphism $\Mr(B) \cong \Mr(B')$ that preserves both meridians and longitudes.

\section{Some questions}
\label{questions}

We have many questions about the ideas discussed in this paper. Here are several.

1. For a classical link, the reduced Alexander module $\Mr(L)$ is a simplification of a module over the ring of Laurent polynomials in $\mu$ variables. For a virtual link, there are several different multivariate Alexander modules that simplify to $\Mr(L)$; for instance see \cite{B, SW}. How do the ideas we have discussed extend to the multivariate setting? In particular, can the formula of Definition \ref{longdef} be modified to define indeterminate longitudes in these modules?

2. In this paper we have focused on modules and their elements, but the ideas we have discussed were developed as we studied the connection between $\Mr(L)$ and the medial quandle of $L$ \cite{mvaq2, mvaq4}. How is the medial quandle connected to $\Me(L)$? In particular, does $\Me(L)$ determine the medial quandle?

3. How do answers to 1.\ and 2.\ extend to the multivariate Alexander quandles of \cite{mvaq3,mvaq1}?

4. Considering that Milnor's $\bar \mu$-invariants are defined using longitudes \cite{Mi}, it is natural to guess that there is a connection with $\Me$. This guess is supported by the fact that the two versions of the Borromean rings are distinguished both by $\Me$ and by $\bar \mu (1,2,3)$ \cite{JLWW}. What is the precise connection between $\Me$ and the $\bar \mu$-invariants? Does the connection also hold for the extended $\bar \mu$-invariants of virtual links recently introduced by Chrisman \cite{Ch}?

5. Lemma \ref{linklem} implies that every classical link has
\[
\epsilon^\mu \varphi_L \left( \sum_{i=1}^\mu \chi_i(L) \right) = 0.
\]
In fact, every classical link we have analyzed has $\sum_{i=1}^\mu \chi_i(L) = 0$. Does the latter equality hold for all classical links?

\section*{Acknowledgment}
The paper was improved by the advice of an anonymous reader.

\end{document}